\documentclass[twoside,openright,11pt,a4paper,epsf]{article}
\usepackage{latexsym}
\usepackage[T1]{fontenc}
\usepackage{amsmath}
\usepackage{amsthm}
\usepackage[french]{babel}
\usepackage{amsfonts} 
\usepackage{amssymb}
\usepackage{vmargin}
\usepackage{fancyheadings}
\usepackage{epsfig}
\usepackage[all]{xy}

\newcommand{\color}[6]{}

\newcommand{\R}{\mathbb{R}}
\newcommand{\N}{\mathbb{N}}

\newcommand{\C}{\mathbb{C}}
\newcommand{\D}{\mathbb{D}}

\newcommand{\cp}{\mathbb{\mathcal P}}
\newcommand{\rc}{\mathbb{\mathcal R}}
\newcommand{\cc}{\mathbb{\mathcal C}}
\newcommand{\bc}{\mathcal B}
\newcommand{\ec}{\mathcal E}

\newcommand{\SC}{\mathcal{SC}}
\newcommand{\spc}{\mathcal{SPC}}

\newcommand{\re}{\text{Re}\,}

\newcommand{\grad}{\vec{\nabla}}

\newcommand{\priv}{\backslash}
\newcommand{\lra}{\longrightarrow}

\newcommand{\om}{\Omega}
\newcommand{\eps}{\varepsilon}
\newcommand{\F}{\mathcal{F}}
\newcommand{\K}{\mathcal{K}}
\renewcommand{\L}{\mathcal{L}}
\newcommand{\psh}{$p.s.h$ }
\newcommand{\jac}{\text{Jac}\,}
\newcommand{\CR}{\text{\tiny CR} }
\newcommand{\wdt}[1]{\widetilde{#1}}
\newcommand{\cqfd}{\hfill $\square$ \vspace{0.1cm}\\ }
\newcommand{\sbull}{{\tiny $\bullet$}}
\newcommand{\ds}{\displaystyle}

\newtheorem*{thm}{Théorème}
\newtheorem*{thm1}{Théorème 1}
\newtheorem*{thm2}{Théorème 1'}  
\newtheorem{definition}{D\'efinition}[section]

\newtheorem{lemme}[definition]{Lemme}
\newtheorem{theoreme}[definition]{Th\'eor\`eme}
\newtheorem{corollaire}[definition]{Corollaire}
\newtheorem{prop}[definition]{Proposition}
\newtheorem{rque}[definition]{Remarque}

\newcounter{faitmoins}

\title{Contractibilité et sphéricité des hypersurfaces strictement 
pseudoconvexes.}
\author{Emmanuel Opshtein}
\begin{document}
\maketitle
\begin{abstract}
Nous proposons une preuve élémentaire de la caractérisation des hypersurfaces 
sphériques comme les seules hypersurfaces strictement pseudoconvexes de $\C^k$ 
possèdant des germes de difféomorphismes contractants. A la différence de
démonstrations antérieures, notre approche n'utilise pas la théorie de 
Chern-Moser et repose entièrement sur les techniques de dilatations de 
coordonnées de Pinchuk. Nous donnons également quelques applications de ce
résultat.
\end{abstract}

\section*{Introduction.}
Dans cet article, nous étudions la géométrie $\text{CR}$ des hypersurfaces
strictement pseudoconvexes de $\C^k$ au moyen des techniques de
dilatations de coordonnées introduites par Pinchuk (\cite{pinchuk4}, 1978). Nous
nous intéressons à la question suivante : un germe d'hypersurface strictement
pseudoconvexe de $\C^k$ peut-il être contracté  sur l'un de ses points par une
suite de difféomorphismes $\text{CR}$ qui préserve l'hypersurface ? 

Les travaux de Chern et Moser (\cite{chemos}, 1974) conduisent à une liste 
complète d'invariants $\text{CR}$ des hypersurfaces Levi non-dégénérées.
Dans une situation contractante ces invariants sont constants et cela suggère 
une certaine homogénéité de la surface. Ainsi, Webster (\cite{webster},
1974), Burns-Shnider (\cite{burshn}, 1977),  Beloshapka (\cite{beloshapka},
1977) et Loboda (\cite{loboda}, 1979) établirent le résultat suivant :
\begin{thm1} Soient $S$ et  $S'$ deux hypersurfaces
  strictement  pseudoconvexes de $\C^k$ (éventuellement à bord).
  S'il existe une suite $(f_n)_n$ de difféomorphismes $\text{CR}$
  de $S$ à valeurs dans $S'$ qui converge localement
  uniformément vers un point de $S'$ alors $S$ est localement sphérique,
 \textit{i.e.} localement  $\text{CR}$-difféomorphe à un ouvert de la sphère 
euclidienne.
\end{thm1}
La démarche de Webster et Burns-Shnider est la suivante. La résolution du
problème de $G$-structure effectuée dans \cite{chemos} produit un tenseur
$\text{CR}$-invariant nul si et seulement si l'hypersurface est
sphérique. Dans le cas non-sphérique, une normalisation de cet invariant
conduit à la réduction du groupe structural de l'hypersurface à
$U(k-1)$. Plus concrètement, ceci montre l'existence d'une distance
$\text{CR}$-invariante sur les hypersurfaces non-sphériques et interdit donc
tout comportement dynamiquement contractant. Beloshapka et Loboda déduisent 
quant à eux ce théorème d'une étude de l'action des difféomorphismes $\text{CR}$
sur les chaînes, et plus précisément d'un résultat de linéarisation locale du
groupe d'automorphismes $\text{CR}$ de l'hypersurface ainsi qu'au théorème 1.  

Un analogue de ce théorème pour les biholomorphismes des domaines est le
théorème de Wong-Rosay. En 1977, Wong \cite{wong} caractérisa la
boule unité parmi les domaines strictement pseudoconvexes bornés de $\C^k$
comme étant le seul à posséder un groupe d'automorphismes non-compact. Rosay
donna une version semi-locale
de l'argument de Wong (\cite{rosay}, 1979), puis Pinchuk en donna un
démonstration basée sur les techniques de dilatations de coordonnées.
\begin{thm}[Wong-Rosay]
Soit $\om$ un domaine  de $\C^k$  et $(f_n)_n$ une suite d'automorphismes 
de $\om$. S'il existe un point $z_0\in \om$ dont l'orbite $\{f_n(z_0)\}_n$ 
accumule un point au voisinage duquel $b\om$ est lisse et strictement
pseudoconvexe alors $\om$ est biholomorphiquement équivalent à la
boule unité de $\C^k$. 
\end{thm}
Les techniques de dilatations ont également permis des généralisations
 de ce résultat au cas faiblement pseudoconvexe. Un théorème analogue 
a été obtenu par  Berteloot \cite{berteloot2} lorsque le point d'impact 
est de type fini dans $\C^2$  et par  Gaussier \cite{gaussier} dans le cas 
convexe de type fini en toutes dimensions. Bedford-Pinchuk avaient auparavant 
caractérisé les domaines dont la frontière est globalement de type fini dans $\C^2$ 
ou convexe de type fini dans $\C^k$ et dont le groupe d'automorphisme est non-compact
(\cite{bedpin2,bedpin1,bedpin3}). 
Ainsi notre étude est susceptible de s'étendre au cas faiblement pseudoconvexe.

Résumons Le principe de notre preuve. Après les avoir 
prolongées en des applications holomorphes entre des domaines soutenus par 
$S$ et $S'$,
nous renormalisons la suite de ces applications en les composant au but par
des dilatations anisotropes bien choisies (partie 1). Le problème est alors de
voir que la suite renormalisée converge vers un biholomorphisme entre le domaine 
source et un ouvert de la boule soutenu par une portion de la sphère (partie 3). 
La difficulté principale est de s'assurer de la non-dégénérescence de 
l'application limite. Ce problème ne se présente pas dans la situation de 
Wong-Rosay. Nous le surmontons en quantifiant
la dilatation des applications $\text{CR}$ (partie 2). Ceci se révèlera 
également utile pour les applications que nous 
présentons dans la partie 4.

Je tiens à remercier François Berteloot pour son aide lors de la rédaction 
de cet article.
\paragraph{Notations et conventions :}  
\begin{itemize}
\item[\sbull] Pour $z\in \C^k$ ($k\geq 2$), nous noterons $z=(z_1,z')$, où $z_1\in
  \C$ et $z'\in \C^{k-1}$.
\item[\sbull] Pour une hypersurface réelle lisse $M$ de $\C^k$, on note
  $\spc(M)$ l'ensemble des points de stricte pseudoconvexité de $M$.
\item[\sbull] Lorsque $M$ est une hypersurface réelle lisse
strictement convexe de $\C^k$ et $q$ est un point de $M$, on note :
\begin{itemize} 
\item $\vec N(q)$ le vecteur unitaire normal à $M$ en $q$ dirigé vers la
  partie convexe de $\C^k\priv M$, 
\item $\Lambda_q$ la forme linéaire $\Lambda_q(\cdot)=\langle \vec
  N(q),\cdot\rangle$,   
\item $q_\eps$ le point $q+\eps\vec N(q)$,
\item $B_\eps(q)$ la boule centrée en $q_\eps$ et de rayon $\eps$; cette
  boule est tangente à $M$ en $q$,
\item $B^+_\eps(q):=B_\eps(q)\cap \{\Lambda_q\geq \eps\}$. 
\end{itemize} 
\item[\sbull] Lorsque $\om$ et $\om'$ sont des domaines de $\C^k$ et $U$
  est inclus dans $b\om$, nous dirons qu'une application $F$ de $\om$ dans 
$\om'$ est propre relativement à $U$ si $F\in \cc^\infty(\om\cup U)$ et 
$F(U)\subset b\om'$.
\item[\sbull] Enfin, lorsque $D$ est un domaine de $\C^k$, et 
$\lambda\in\R$, nous noterons 
$$
\begin{array}{ll}
D_\lambda:=D\cap \{\re z_1<\lambda\} & D_{\lambda}^+:=D\cap \{\re z_1\geq
\lambda\}\\
 {[bD]}_\lambda:=bD\cap \{\re z_1<\lambda\} & {[bD]}_\lambda^+:=bD\cap
\{\re z_1\geq\lambda\} .\\
\end{array}
$$
\end{itemize}
\section{Renormalisation de la suite $f_n$.}
Commençons par décrire et simplifier la situation géométrique sous-jacente
aux hypothèses du théorème $1$.

Appelons $a$ le point de $S'$ vers lequel $(f_n)_n$ converge. Soit $p$ un
point de $S$. Le but est de trouver un voisinage de $p$ dans $S$ qui est 
$\text{CR}$-difféomorphe à un ouvert de la sphère. Ce problème ne dépend évidemment
pas des systèmes de coordonnées holomorphes en $p$ et en $a$. Comme les
surfaces $S$ et $S'$ sont strictement pseudoconvexes, on peut les
supposer strictement {\it convexes} sur des voisinages $U_2$ et $U'$ de
$p$ et $a$ respectivement. De même, on peut supposer $\Lambda_p=\Lambda_a=\re
z_1$, $U_2=S\cap \{\re z_1< 2\}$ et $U'=S'\cap \{\re z_1< 1\}$. Il
existe alors deux domaines bornés de $\C^k$ strictement convexes, $\om$ et
$\om'$ tels que $U_2= [b\om]_2$ et $U'= [b\om']_1$. 

Quitte à éliminer de la suite $(f_n)_n$ un nombre fini d'éléments, on peut
supposer que $f_n([b\om]_2)\subset U'$ pour tout $n\in \N$. 
Les théorèmes classiques de prolongement des applications $\text{CR}$
(voir \cite{bogess}, chap. 15) montrent que $f_n$ se prolonge en une application
holomorphe de $\om_2$ dans $\C^k$ que nous 
noterons abusivement $f_n$. Comme $\om_1'$ possède une fonction \psh\
définissante $\rho$, le principe du maximum appliqué 
aux restrictions de $\rho\circ f_n$ sur les hyperplans complexes $\{z_1=c\}$
assure que $f_n(\om_2)\subset \om_1'$. Un argument similaire prouve que la
suite $(f_n)$ converge uniformément vers $a$ sur $\om_2$.

\begin{figure}[h!]
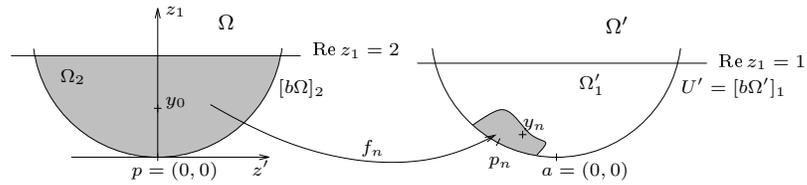

\begin{center}
\input sit1.pstex_t
\end{center}
\caption{Situation géométrique du théorème $1$.}
\label{situation1}
\end{figure}

Posons enfin $y_0:=(1,0)$, $y_n:=f_n(y_0)$, $p_n:=f_n(p)$, de sorte que $y_n$
et $p_n$ tendent vers $a$. La figure \ref{situation1} illustre la situation.
\subsection{Dilatations de coordonnées.}
\begin{figure}[h]
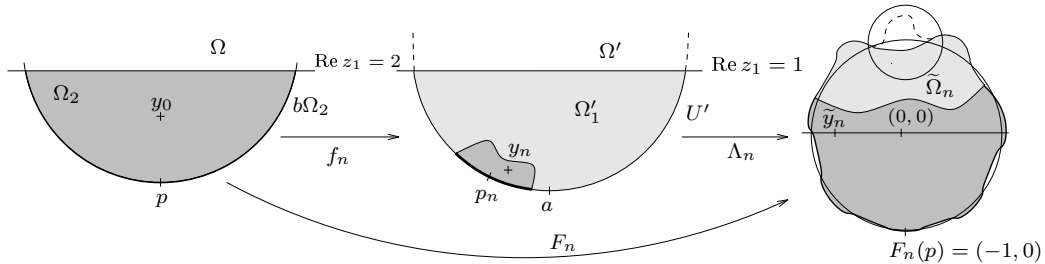

\begin{center}
\input dilat.pstex_t
\end{center}
\vspace{-0.1cm}
\caption{Situation après remise à l'échelle ``au temps'' $n$.}
\label{situation2}
\end{figure}

Nous décrivons une technique de dilatation de coordonnées naturellement
associée au couple $(p_n,y_n)$. Elle procède d'une légère modification de
la méthode introduite par Pinchuk (voir \cite{pinchuk,berteloot}). 
On peut supposer que $U'$ 
possède une fonction définissante $\chi$ donnée par : 
$$
\chi(z):=-\re z_1+\Vert z\Vert^2+a(z)\Vert z\Vert^3 \text{ où } a(z)\in \mathcal
O(1).  
$$ 
L'opération de remise à l'échelle en $q\in U'$ de paramètre $\eps$, notée 
$\rc_{q,\eps}$, consiste en
une composition de deux biholomorphismes de $\C^k$, $\Psi_q$ et $D_\eps$, et
d'une transformation  birationnelle $\Phi$ de $\C^k$ :
$$
\rc_{q,\eps}=\Phi\circ D_\eps\circ \Psi_q.
$$
\sbull\ $\Psi_q$ est la composition d'une translation qui envoie $q$ sur
l'origine et d'une transformation unitaire de $\C^k$ qui rend horizontal le
plan tangent à $b\om$ en $q$ (\textit{i.e.} défini par $\re z_1=0$). Dans ces
coordonnées, $\chi_q:=\chi\circ \Psi_q^{-1}$ est de la forme :
\begin{equation}\label{eqdil1}
\chi_q(z)=-\re z_1+Q_q(z)+a_q(z)\Vert z\Vert^3  
\end{equation}
où $Q_q$ est le hessien réel de $\chi_q$ en $0$. Il s'agit d'une forme
quadratique réelle définie positive sur $\C^k=\R^{2k}$ 
car $U'$ est strictement convexe. La fonction $a_q(z)$ est majorée sur
$\om_1'$ indépendamment de $q\in U'$ car $b\om'$ est $\cc^3$. Remarquons que
$\Psi_a=Id$, $Q_a=\Vert \cdot\Vert$ et que les 
applications $q\mapsto \Psi_q$ et $q\mapsto Q_q$ sont continues en raison du
caractère lisse de $U'$.\vspace{0.2cm}\\    
\sbull\ $D_\eps$ est la dilatation anisotrope $\C$-linéaire 
$$
D_\eps:(z_1,z')\longmapsto \left(\frac{z_1}{\eps},\frac{z'}{\sqrt \eps}\right).
$$
Il s'agit d'un changement d'échelle naturel associé à la stricte
pseudoconvexité. Un calcul simple montre que la fonction définissante
$\chi_{q,\eps}$ de $D_\eps\circ \Psi_q(U')$ définie par
$\ds\chi_{q,\eps}:=\frac{1}{\eps} \chi_q\circ D_\eps^{-1}$ s'écrit :
$$
\chi_{q,\eps}(z_1,z') = -\re z_1+Q_q\big((0,z')\big) +O(\sqrt \eps).
$$  
\sbull\ $\Phi$ est la transformation de Cayley 
$$
\Phi:(z_1,z')\longmapsto\left(\frac{z_1-1}{z_1+1},\frac{2z'}{z_1+1}\right)
$$ 
qui transforme $\Sigma:=\{\re z_1\geq \Vert z'\Vert^2\}$ en $B:=\{| z_1|^2+\Vert
z'\Vert^2\leq 1\}$. Remarquons que $D_\eps\circ\Psi_q(\om')$ est inclus dans le
demi-espace $\{\re z_1>0\}$ en raison de la convexité de $\om'$. La transformation 
$\rc_{q,\eps}=\Phi\circ D_\eps\circ \Psi_q$ est par conséquent bien définie
sur $\om'$. 
\paragraph{}Nous allons maintenant définir une suite de remises à l'échelle 
$\rc_n:= \rc_{p_n,\eps_n}$ associée à la suite $(p_n,y_n)$ ; voir la figure
 \ref{situation2}. A cet effet, observons que  
$\{\rc_{q,\eps}^{-1}(\{\re z_1=0\}),\eps>0\}$
définit un feuilletage de $\C^k\priv T_q^\C b\om'$ dont les feuilles 
sont des cylindres euclidiens centrés en $T_q^\C b\om'$ :
$$
F_{q,\eps}:=\rc_{q,\eps}^{-1}(\{\re z_1=0\})=\{z\in \C^k,\; d(z,T_q^\C
b\om')=\eps\}.  
$$
En raison de la convexité de $\om'$, tout point $z\in\om'$ appartient
à une unique feuille $F_{q,\eps_q(z)}$, $\eps_q(z)>0$. Ceci nous permet
d'énoncer la définition suivante :
\begin{definition}
La remise à l'échelle $\rc_n$ associée à la suite $(p_n,y_n)$ est définie 
par $\rc_n:=\rc_{p_n,\eps_n}$ où $\eps_n:=\eps_{p_n}(y_n)$ est l'unique réel 
tel que $y_n\in F_{p_n,\eps_n}$. On pose $F_n:=\rc_n\circ f_n$, 
$\wdt\om_n:=\rc_n(\om')$ et $\wdt y_n:=\rc_n(y_n)$.   
\end{definition}
Grâce à ces choix nous disposons d'une suite d'applications $F_n$ telles que 
(voir fig. \ref{situation2}) :
$$
\begin{array}{rcccl}
F_n & : & \om_2 &\lra & \wdt\om_n \\
& & p & \longmapsto & a=(-1,0)\\
& & y_0 & \longmapsto & \wdt y_n \in \{\re z_1=0\}.\\
\end{array}
$$

Notons que la convergence de $p_n$ et $y_n$ vers $a$ implique que 
$\eps_{p_n}(y_n)$ tend vers $0$. Il est bien connu que la 
suite de domaines $\wdt \om_n$ tend vers la boule $B$ au sens $\cc^2$ sur les
compacts de $\C^{k+1}\priv \{z_1=1\}$. L'objet de la proposition suivante
 est de mieux contrôler cette convergence. Ceci se révèlera
utile dans la troisième partie. Commençons par observer que $\om'$ est 
inclus dans toutes les paraboles tangentes à $b\om'$ de courbure assez petite :
\begin{lemme}\label{geodiff}
Il existe une constante $c>0$ telle que $\Psi_q(\om')\subset \{\re
z_1>c\Vert z\Vert^2\}$ pour tout $q\in b\om'$.
\end{lemme}
\noindent\underline{Preuve :} Raisonnons par l'absurde,
supposons qu'il existe une suite de réels $c_n$ qui tend vers $0$ et des
points distincts $p_n,p_n'\in b\om'$ tels que $p_n'\in
\cp_n:=\Psi_{p_n}^{-1}\big(\{\re z_1>c_n\Vert z\Vert^2\}\big)$ ($\cp_n$ est
une parabole tangente à $b\om'$ en $p_n$). Comme $b\om'$ est compact, on peut
supposer que $p_n$ et $p_n'$ convergent vers deux points $p,p'\in b\om'$.
 Nous allons montrer que $p\neq p'$. Ceci contredira la stricte convexité de
$\om'$ puisqu'alors $p\neq p'\in T_pb\om'\cap b\om'$ (la suite de paraboles
$\cp_n$ converge sur les compacts vers $T_pb\om'$ car $c_n$ tend vers $0$ et
$p_n$ tend vers $p$).

On peut évidemment
supposer que $p=0$. Rappelons que (\ref{eqdil1}) donne l'expression d'une
fonction définissante de $\Psi_q(U')$ pour $q\in U'$. Pour $q$
suffisamment proche de l'origine, on peut supposer  que $Q_q(z)\geq 1/2\Vert
z\Vert^2$. Il existe donc un voisinage $B$ de l'origine dans $\C^k$ tel que 
$$
\Psi_q(\om'\cap B)\subset \{\re z_1>\frac{\Vert z\Vert^2}{3}\}.
$$
Comme $c_n$ tend vers $0$, le point $p'_n$ est hors de $B$ donc $p'\neq
p$. \cqfd 
\indent Nous dirons qu'une suite de domaines $D_n$ converge vers un
domaine $D$ au sens $\mathcal C^2$ sur un compact $K$ si $K\cap
\overline{D_n}$ tend vers $K\cap \overline{D}$ au sens de la topologie de
Hausdorff et si les fonctions définissantes $\chi_{i,n}$ de $D_n$ sur les
ouverts $U_i$ d'un recouvrement de $D$ convergent en norme $\cc^2$ vers des
fonctions définissantes $\chi_i$ de $D\cap U_i$.  
%\begin{itemize}
%\item $K\cap \overline{D_n}$ tend vers $K\cap \overline{D}$ au sens de la
%  topologie de Hausdorff.
%\item Lorsque $bD$ est recouvert par des ouverts $U_i$ de $\C^k$ sur lesquels il
%  possède des fonctions définissantes $\chi_i$ $\mathcal C^2$, alors  $bD_n$ a
%  des fonctions définissantes $\chi_{i,n}$  $\mathcal C^2$-proches de $\chi_i$ 
%  sur  $U_i$ pour $n$ suffisamment grand. 
%\end{itemize}
\begin{prop}\label{cvcool}
La suite de domaines $(\wdt\om_n)$ converge vers $B$ au sens $\mathcal C^2$ 
sur les compacts de $\C^k\priv\{(1,0)\}$.
\end{prop}
\noindent\underline{Preuve :} Les domaines $\wdt \om_n$ sont de la forme 
$$
\wdt \om_n=\Phi\circ D_{\eps_n}\circ \Psi_{p_n}(\om')
$$
où $\eps_n$ tend vers $0$. On sait que la suite
de domaines $D_{\eps_n}\circ \Psi_{p_n}(\om')$ tend vers $\Sigma$ au sens
$\mathcal C^2$ sur les compacts de $\C^k$ (voir \cite{pinchuk}) et que la
transformation de Cayley $\Phi$ transforme $\Sigma$ en $B$ en faisant 
correspondre l'infini à $\{z_1=1\}$. Il suffit donc de vérifier qu'il existe 
une fonction $r(\delta)$ tendant vers $0$ avec $\delta$ telle que :
$$
\wdt \om_n\cap \{|z_1-1|<\delta\}\subset B[(1,0),r(\delta)] \hspace{0.5cm} 
\forall n\in \N.
$$
Or ceci résulte du lemme \ref{geodiff}. En effet, d'après ce lemme, il existe 
$c>0$ tel que $\Psi_{p_n}(\om')\subset \{\re z_1>c\Vert z'\Vert^2\}$ pour tout 
$n\in \N$. En outre, on vérifie aisément que 
$\Phi\circ D_{\eps_n}(\{\re z_1>c\Vert z'\Vert^2\})=E_c:=\{|z_1|^2+c\Vert 
z'\Vert^2<1\}$. L'inclusion attendue résulte alors immédiatement de
$$
\wdt \om_n\cap \{|z_1-1|<\delta\}\subset E_c\cap \{|z_1-1|<\delta\}.
$$\hfill $\square$\vspace{-0.1cm}

La figure \ref{situation2} décrit le procédé de remise à l'échelle. 

 Observons que la normalité de la famille $(F_n)$ résulte 
de la proposition \ref{cvcool} puisque celle-ci
implique que $\wdt\om_n$ est une suite de domaines uniformément bornés.
\begin{corollaire}\label{normal}
La famille $(F_n)$ est normale. Quitte à extraire, on peut supposer que $F_n$
converge uniformément vers $F:\om_2\lra \overline{B}$ sur les compacts de $\om_2$. 
De plus l'alternative suivante a lieu : soit $F$ est à valeurs dans $B$ soit 
$F(\om_2)$ est un point de $bB$.  
\end{corollaire}

\begin{rque}\label{bord}
L'image de $U'$ par la remise à l'échelle en $q\simeq a$ et de paramètre
$\eps\ll 1$ est une hypersurface à bord de $\C^k$. Son bord est non vide et
inclus dans un voisinage de $(1,0)$ dont le diamètre tend vers $0$ avec
$\eps$.  
\end{rque}

\subsection{Réduction du problème.}
La preuve du théorème 1 consiste à montrer que $F_n$ converge vers un
biholomorphisme d'un voisinage de $p$ dans $\overline{\om}$ sur un voisinage
de $(-1,0)$ dans $\overline{B}$. Nous réduisons ici  ce problème à  
l'étude de la convergence des applications $F_n$ et de leurs inverses $G_n$. 
Montrons tout d'abord comment définir ces inverses.
\begin{lemme}\label{remplissage}
Soient $D$ et $D'$ deux domaines à bords lisses de $\C^k$. On suppose que
$D_1$ et $D'_1$ sont connexes, non vides et que $[bD]_1$, $[bD']_1$ sont
strictement  convexes. Soit $F:D_1\lra D'$ une application holomorphe propre 
relativement à $[bD]_1$ telle que $F([bD]_1)\subset \spc(bD')$. Alors $F$ est 
un difféomorphisme local sur $D_1$.  De plus, si $F$ est injective sur $[bD]_1$
et si $F([bD]_1)\supset[bD']_1$ alors $F(D_1)\supset D'_1$ et il existe 
un inverse à droite de $F$ défini sur $D'_1$. 
\end{lemme}
\noindent\underline{Preuve :} Notons tout d'abord que si $F$ n'était pas un 
difféomorphisme local sur $D_1$ alors le lieu critique $\{\jac F=0\}$ serait un 
ensemble analytique dont l'intersection avec $[bD]_1\subset \spc(bD)$ serait 
vide \cite{pinchuk3}. La fonction $-\re z_1$ violerait le principe 
du maximum sur cet ensemble analytique.

Comme $[bD]_1$ est formé de points de stricte
convexité et $F$ y est injective, $F_{|[bD]_1}$ est un difféomorphisme $\text{CR}$
sur son image, qui contient $[bD']_1$. Il existe donc un difféomorphisme $\text{CR}$ 
$G$ de $[bD']_1$ dans $[bD]_1$ tel que $F\circ G=Id$. Celui-ci se prolonge en une 
application $\overline{G}:D_1'\lra D$ holomorphe propre relativement à $[bD']_1$ 
dont l'image est incluse dans $D_1$ d'après le principe du maximum. L'application 
$F\circ \overline{G}$ est donc bien définie sur $D'_1\cup [bD']_1$, et coïncide 
avec l'identité sur $[bD']_1$. Le principe d'unicité montre que  
$F\circ \overline{G}=Id_{D'_1}$. On en conclut que $F(D_1)\supset D'_1$ et que 
$\overline{G}$ est l'inverse de $F$ sur $D_1'$.\cqfd
Ce lemme assure l'existence d'un inverse $G_n$ à $F_n$ naturellement défini
sur un ensemble maximal du type $\wdt \om_{n,t}$. La définition suivante met en 
relief les propriétés de ces applications que nous aurons à vérifier :
\begin{definition}\label{ko}
Nous désignerons par $\K(F_n)$, $O(F_n)$, $\K(G_n)$, $O(G_n)$ les propriétés 
suivantes :
\begin{description}
\item[\sbull] $\K(F_n)$ : $F_n$ converge vers une application $F$ et 
$F(\om_2)\subset B$,
\item[\sbull] $O(F_n)$ : il existe un réel $\lambda$ strictement positif  tel que
  $F_n(\om_1)\supset \wdt\om_{n,\lambda}$ pour tout  $n\in \N$,
\item[\sbull] $\K(G_n)$ : $G_n$ converge vers une application $G$ et 
$G(B_\lambda)\subset \om_2$,
\item[\sbull] $O(G_n)$ : il existe une fonction croissante continue 
$\lambda':]0,\lambda]\to\R^+_*$  telle que \\
$G_n(\wdt\om_{n,t})\supset \om_{\lambda'(t)}$ pour tout $n$.  
\end{description}
\end{definition} 
Nous sommes maintenant en mesure de formuler la réduction annoncée :   
\begin{prop}\label{reduction}
Si la suite de renormalisation $(F_n)_n$ et la suite des inverses $(G_n)_n$ 
vérifient les propriétés $\K(F_n)$, $O(F_n)$, 
$\K(G_n)$, $O(G_n)$ alors l'application limite $F:\om_2\to B$ se prolonge en un 
$CR$-difféomorphisme entre un voisinage de $p$ dans $[b\om]_1$  et un voisinage 
de $(-1,0)$ dans $bB$. 
\end{prop}
\noindent\underline{Preuve :}
D'après l'hypothèse $O(F_n)$ et le lemme \ref{remplissage}, $G_n$ est définie 
sur $\wdt \om_{n,\lambda}$ et y vérifie :
\begin{equation}\label{reduc1}
F_n\circ G_n=Id.
\end{equation} 
On déduit immédiatement de (\ref{reduc1}) et de l'hypothèse $O(G_n)$ que 
$F_n(\om_{\lambda'(t')})\subset \wdt \om_{n,t'}$ pour $t'<t\leq \lambda$. 
Compte tenu de l'hypothèse $\K(F_n)$, cela donne à la limite : 
$F(\om_{\lambda'(t')})\subset \overline{B}_{t'}\cap B\subset B_t$. On en déduit 
par continuité de $\lambda'$ que :
\begin{equation}\label{reduc2}
F\big(\om_{\lambda'(t)}\big)\subset B_t \hspace{0.3cm}\text{ pour } t\leq \lambda.
\end{equation}
Il résulte également de (\ref{reduc1}) que $G_n\circ F_n=Id$ sur $G_n\big(\wdt 
\om_{n,\lambda}\big)$ et donc sur $\om_{\lambda'(\lambda)}$ en vertu de l'hypothèse 
$O(G_n)$ :
\begin{equation}\label{reduc3}
G_n\circ F_n=Id \hspace{0.3cm}
\text{ sur } \om_{\lambda'(t)} \text{ pour } t\leq \lambda.
\end{equation}
L'hypothèse $\K(G_n)$ permet de passer à la limite dans (\ref{reduc1}) car 
$F_n$ converge uniformément vers $F$ sur les compacts de $\om_2$. Il vient alors :
\begin{equation}\label{reduc4}
F\circ G=Id \hspace{0.3cm}\text{ sur } B_t\text{ pour } t\leq \lambda.
\end{equation}
L'inclusion (\ref{reduc2}) permet de passer à la limite dans (\ref{reduc3}) :
\begin{equation}\label{reduc5}
G\circ F=Id \text{ sur } \om_{\lambda'(t)} \text{ pour } t\leq \lambda.
\end{equation}
L'identité (\ref{reduc4}) montre que $F$ induit un biholomorphisme entre 
$G(B_\lambda)$ et $B_\lambda$ :
$$
F:G(B_\lambda)\overset{\sim }{\lra} B_\lambda.
$$
En outre, l'inclusion (\ref{reduc2}) et l'identité (\ref{reduc5}) montrent 
que $G(B_\lambda)\supset \om_{\lambda'(\lambda)}$. On est alors en mesure de
montrer que $F$ se prolonge en un difféomorphisme d'un voisinage de $p$ dans 
$b\om$ sur un voisinage de $(-1,0)$ dans $bB$. Les frontières étant strictement 
pseudoconvexes en $p$ et $(-1,0)$, un théorème de Bell \cite{bell} stipule que 
le prolongement différentiable est une conséquence du prolongement continu. 
Fixons $-1<t\leq \lambda$ arbitrairement proche de $-1$. 
Pour $n$ assez grand $x_n\in \om_{\lambda'(t)}$ et donc, compte tenu de 
(\ref{reduc2}), $F(x_n)\in B_t$. Ainsi $F(x_n)$ tend vers $(-1,0)$ lorsque $x_n$
tend vers $p$ et cela entraîne le prolongement hölder continu de $F$ à $b\om$
autour de $p$ \cite{berteloot3}.\cqfd

Notre objectif est de vérifier les propriétés ``$\K$'' et ``$O$'' pour les 
suites $F_n$ et $G_n$. Il s'agit de propriétés de compacité de suites 
d'applications respectivement à l'intérieur et au bord des domaines. Nous 
mettons au point dans la partie suivante des estimations sur les dérivées et les 
images des applications $\text{CR}$ qui nous permettront de les établir.
\section{Préliminaires techniques sur les applications $\text{CR}$.}\label{prelim}
Nous introduisons dans cette partie une distance qui présente l'avantage
d'être plus adaptée à la géométrie des applications $\text{CR}$ tout en
définissant la même topologie que la distance euclidienne sur les
hypersurfaces strictement pseudoconvexes. 

Soit $X$ une sous-variété lisse de $\C^k$, $x$ et $y$ deux points de $X$. On
appelle chemin complexe entre $x$ et $y$ tout chemin $\gamma :[0,l]\lra X$
joignant $x$ à $y$, $\mathcal C^1$ par morceaux, et tel que
$\dot{\gamma}(t)\in T^\C_{\gamma(t)} S$ là où cela fait sens. On note
$\ell(\gamma)$ la longueur euclidienne d'un tel chemin. 
Pour $x,y\in X$, on définit la distance $\text{CR}$ entre $x$ et $y$ par :
$$
d^{\text{\tiny CR}}_X(x,y):=\text{inf}\{\ell(\gamma),\;\gamma \text{ chemin
  complexe entre  $x$ et $y$}\}. 
$$
En l'absence d'ambiguïté, on ne précisera pas la dépendance en $X$ de
$d^{\text{\tiny CR}}$. Les boules correspondantes de centre $x$ et de
rayon $\delta$ sont notées $B_X^{\text{\tiny CR}}(x,\delta)$. 

Les distances euclidienne et
$\text{CR}$ sont fortement reliées dans les hypersurfaces minimales
(\textit{i.e.} ne contenant pas  de bout d'hypersurface holomorphe). En
particulier, les grandes boules $\text{CR}$  contiennent de grandes boules
euclidiennes. Pour les hypersurfaces strictement pseudoconvexes,
cette propriété est un cas particulier très simple du théorème 4 de \cite{nsw}.
\begin{prop}\label{minimale}
La topologie associée à la distance $\text{CR}$ coïncide avec la
topologie usuelle sur les hypersurfaces strictement pseudoconvexes lisses de
$\C^k$. En particulier, tout ouvert borné d'une telle hypersurface est
$d^{\text{\tiny CR}}$-borné.  
\end{prop}
Le résultat principal de cette partie est une version quantifiée du principe
suivant :  
\begin{center}
"Une application $\text{CR}$ est une dilatation pour la distance $\text{CR}$."
\end{center}
Compte tenu de la proposition \ref{minimale}, ce principe implique que les
images des applications $\text{CR}$ contiennent naturellement de "grands"
ouverts euclidiens. Il est régi par le lemme de Hopf, qui permet de
minorer les dérivées des applications holomorphes propres. 
\begin{lemme}[Hopf] \label{hopf}
Soit $D$ un domaine de $\C^k$ et $\chi\in \cc^1(\overline{D})$ une fonction
\psh\ négative sur $D$ et nulle en $q\in bD$. Soit $\eps$ tel que
$B_\eps(q)\subset D$ et $M:=\min\{|\chi(z)|,z\in
B_\eps^+(q)\}$. Alors $\ds\chi(q_t)\leq -\frac{Mt}{8\eps}$ pour
$t\leq\eps$ et donc $\ds\Vert\grad \chi(q)\Vert\geq \frac{M}{8\eps}$.
\end{lemme} 
\noindent\underline{Preuve :} Soit $\D_\eps$ le paramétrage unitaire du disque
complexe de centre $q_\eps$ passant par $q$, qui envoie $q_\eps$ sur $0$ et
$q$ sur $-\eps$. Dans ce paramétrage,  $B^+_\eps(q)\cap \D_\eps$
correspond à $\D_\eps\cap \{\re z>0\}$. La restriction $\wdt \chi$ de $\chi$
à $\D_\eps$ est une fonction sous-harmonique négative sur $\D_\eps$ vérifiant
$\chi(q_t)=\wdt\chi(-\eps+t)$, $\wdt \chi(\eps e^{i\theta})\leq -M$ pour
$\theta\in[-\frac{\pi}{2},\frac{\pi}{2}]$ et $\wdt\chi(-\eps)=0$. Un calcul facile
sur le noyau de Poisson donne donc : 
$$
\begin{array}{rl} 
\chi(q_t)=\wdt\chi(-\eps+t)\leq &\ds\frac{1}{2\pi}\int_{-\pi}^\pi
\frac{\eps^2-(\eps-t)^2}{|-\eps+t-\eps e^{i\theta}|^2} \wdt\chi(\eps
e^{i\theta})d\theta \\ 
  \leq &  \ds -\frac{M}{2\pi}\int_{-\frac{\pi}{2}}^\frac{\pi}{2}
\frac{\eps^2-(\eps-t)^2}{|-\eps+t-\eps e^{i\theta}|^2}\, d\theta 
 \leq \ds -\frac{M}{2}\frac{t(2\eps-t)}{4\eps^2} 
 \leq  \ds -\frac{M t}{8\eps}.
\end{array}
$$
\cqfd
Nous combinons à présent le lemme de Hopf à une propriété de transfert d'estimations
des dérivées normales aux dérivées tangentielles complexes pour obtenir la version 
du principe énoncé ci-dessus qui nous sera utile.
\begin{lemme}\label{recouvre}
Soient $\om$, $\wdt\om$ des domaines de $\C^k$ à bords lisses et $U$ un ouvert 
relativement compact dans $\spc(b\om)$. Soit $W$ un ouvert de $\spc(b\wdt \om)$ 
tel que $\Lambda_x$ est positive sur $\wdt \om$ pour tout $x\in 
W$. On note $\kappa_1$ et $\kappa_2$ la
plus petite et la plus grande des valeurs propres de la forme de Levi sur $U$
et $W$. Soit $p\in U$ et $\eps,\tau$ deux réels strictement positifs tels que
$B^{\text{\tiny CR}}(p,\tau)\Subset U$ et  $B_\eps(q)\subset \om$ pour tout
$q\in B^{\text{\tiny CR}}(p,\tau)$.

Soit $F:\om\lra \wdt \om$ une application holomorphe propre relativement à $U$
telle que $F(p)\in W$. On définit 
$$ 
M:=\min\big\{|\Lambda_x(z)|,\; x\in W, \; z\in F\big(\bigcup_{q\in
  B^{\text{\tiny CR}}(p,\tau)} B_\eps^+(q)\big)\big\}.
$$
Alors 
$$
F\big(B^{\text{\tiny CR}}(p,\tau)\big)\supset B^{\text{\tiny
    CR}}_{W}\left(F(p),\left[\frac{\kappa_1}{8\kappa_2}\frac{M}{\eps}
  \right]^{\frac{1}{2}}\tau\right).      
$$
\end{lemme}\vspace{0.3cm}
La preuve comporte deux étapes. Dans la première, le lemme de Hopf
est utilisé pour minorer les dérivées tangentielles complexes de $F$ sur la
préimage de $W$. La seconde étape consiste à intégrer ces estimées.\\
\underline{Preuve :}
{\it Estimations des dérivées de $F$.} Soit $q$ un point de $U$ tel
que $F(q)\in W$. Montrons :
\begin{equation}\label{hopf+}
\Vert F'(q)u\Vert\geq \left[\frac{\kappa_1}{8\kappa_2}\frac{M}{\eps} \right]^{1/2}
\Vert u\Vert,\hspace{0.5cm}\forall u\in T^\C_q b\om.
\end{equation}
Le lemme \ref{hopf} appliqué à la fonction $-\Lambda_{F(q)}\circ F$ fournit
l'estimée $\Vert~\grad ~ [~\Lambda_{F(q)}~\circ ~
F](q)~\Vert~\geq~\frac{M}{8\eps}$. Comme $\grad
\Lambda_{F(q)}\big(F(q)\big)=\vec N(F(q))$, on a :
\begin{equation}\label{eqevident}
n_q(F):=\langle F'(q)\vec N(q),\vec N(F(q))\rangle=
d\Lambda_{F(q)}[ F'(q)\vec N(q)]=\Vert\grad [\Lambda_{F(q)}\circ
F](q)\Vert\geq\frac{M}{8\eps}. 
\end{equation}
Notons $\L(\phi,a,u)$ la forme de Levi en $a$ d'une fonction $\phi$,
appliquée à un vecteur $u$.\\
Comme $W$ est lisse, il existe une fonction définissante
$\psi$ de $\wdt \om$ au voisinage de $F(q)$, qu'on peut supposer de
gradient normé en $F(q)$ : $\grad \psi(F(q))=\vec N(F(q))$. Comme pour 
(\ref{eqevident}), on voit qu'alors $\Vert\grad \psi\circ F(q)\Vert=n_q(F)$. De
plus, $\psi\circ F$ est une fonction définissante en $q$, au voisinage duquel  
$b\om$ est strictement pseudoconvexe. Pour tout vecteur $u\in T^\C_q
b\om$, on a donc compte tenu de (\ref{hopf+}):
\begin{equation}\label{nq}
\L(\psi\circ F,q,u)\geq \kappa_1\Vert \grad \left[\psi\circ
  F\right](q)\Vert\Vert u\Vert^2=\kappa_1n_q(F)\Vert u\Vert^2\geq
  \frac{\kappa_1 M}{8\eps}\Vert u\Vert^2. 
\end{equation}
Par fonctorialité de la forme de Lévi et définition de $\kappa_2$ il vient :
\begin{equation} \label{fonct}
\L(\chi\circ F,q,u)=\L(\chi,F(q),F'(q)u)\leq \kappa_2\Vert F'(q)u\Vert^2.
\end{equation}
Les inégalités (\ref{nq}) et (\ref{fonct}) impliquent (\ref{hopf+}).\vspace{0.2cm}\\   
{\it Conclusion.} Nous intégrons à présent l'inégalité (\ref{hopf+}).
Comme $F(p)$ appartient à $W$, il suffit de vérifier que 
$$
bF\big(B^{\text{\tiny CR}}(p,\tau)\big)\cap B^{\text{\tiny CR}}_{\wdt
  W}\left(F(p),\left[\frac{\kappa_1}{8\kappa_2}\frac{M}{\eps}
  \right]^{1/2}\tau\right)=\emptyset ,  
$$
ou encore que pour tout point $x$ de $bF(B^{\text{\tiny CR}}(p,\tau))\cap W$, 
\begin{equation}\label{eqfin}
d^{\text{\tiny CR}}_{W}(F(p),x)\geq
\left[\frac{\kappa_1}{8\kappa_2}\frac{M}{\eps}\right]^{\frac{1}{2}}\tau.
\end{equation}
Soit $x$ un tel point et $\wdt \gamma$ un chemin complexe de $W$ entre $F(p)$ 
et $x$  paramétré par la longueur d'arc. Comme $F$ est un difféomorphisme
$\text{CR}$ local en tout point de $B^{\text{\tiny CR}}(p,\tau)$, 
on peut relever la composante connexe de $F(p)$ dans $\wdt 
\gamma\cap F(B^{\text{\tiny CR}}(p,\tau))$ en un chemin complexe $\gamma$ tant que 
$\gamma$ ne sort pas de $B^{\text{\tiny CR}}(p,\tau)$.
Précisément, il existe $l\leq \ell(\wdt \gamma)$ et $\gamma:[0,l]\lra 
\overline{B^{\text{\tiny CR}}(p,\tau)}$ joignant $p$ à $bB^{\text{\tiny
    CR}}(p,\tau)$ et tel que $F\circ \gamma(t)=\wdt\gamma(t)$ pour tout
$t\in[0,l]$. 
Comme  $\gamma(t)$ appartient à $U$ et $F(\gamma(t))$ à $W$ pour tout
$t$, l'estimée (\ref{hopf+}) obtenue dans l'étape précédente montre que :
$$
 \Vert
F'(\gamma(t))u\Vert\geq
\left[\frac{\kappa_1}{\kappa_2}\frac{M}{8\eps}\right]^{\frac{1}{2}}\Vert 
u\Vert\hspace{0.5cm}\forall t\in [0,l],\; \forall u\in T^\C_{\gamma(t)}b\om. 
$$
Ainsi,
$$
%\begin{array}{lll}
%\ds d^{\text{\tiny CR}}(F(p),x)+\delta\geq 
\ell(\wdt \gamma) \geq  l  =  \int_0^l\Vert
 \dot{\wdt\gamma}(t)\Vert dt
 = \int_0^l \Vert F'(\gamma(t))\dot{\gamma}(t)\Vert dt \geq 
 \left[\frac{\kappa_1}{\kappa_2}\frac{M}{8\eps}\right]^{\frac{1}{2}}
 \int_0^l\Vert \dot{\gamma}(t)\Vert dt \\ 
  \geq 
 \left[\frac{\kappa_1}{\kappa_2}\frac{M}{8\eps}\right]^{\frac{1}{2}}\ell(\gamma).
%\end{array} 
$$
Puisque $\gamma$ joint $p$ à $bB^{\text{\tiny CR}}(p,\tau)$ et que $\wdt
\gamma$ est un chemin complexe quelconque liant $F(p)$ à $x$, on tire (\ref{eqfin}).
\cqfd

Le lemme \ref{recouvre} permet donc d'estimer la taille de l'image d'une boule 
$B^\CR(p,\tau)$ par une application holomorphe propre en fonction de son 
comportement sur 
$$
A_\eps(p,\tau):=\bigcup_{q\in B^{\text{\tiny CR}}(p,\tau)} B_\eps^+(q).
$$
La proposition suivante détaille les situations dans lesquelles nous allons appliquer 
ce lemme.
\begin{prop}\label{recap}
Soient  $\om$, $\wdt \om$ des domaines de $\C^k$. soient $U$ et $W$ des ouverts 
relativement  compacts dans $\spc(b\om)$ et $\spc(b\wdt \om)$. On suppose que 
$\Lambda_x\geq 0$ sur $\wdt \om$ pour tout $x\in W$.
 
Fixons des points $p\in U$, $a\in W$ et des réels $\eps_0,\tau>0$ tels que 
$B^\CR(p,\tau)\Subset U$ et $A_{\eps_0}(p,\tau)\subset \om$. 

Soit $\F_{p,a}$ la famille des applications holomorphes $F:\om\lra \wdt \om$ propres 
relativement à $U$ vérifiant $F(p)=a$.
\begin{enumerate}
\item A tout compact $K$ de $\wdt \om$ est associée une constante $c(K)>0$  
telle que  :
$$
\forall F\in \F_{p,a} \; :\;\;\big(F(p)\in W\text{ et } F(A_{\eps_0}(p,\tau))\subset K\big)
\Longrightarrow F(B^\CR(p,\tau))\supset B_W^\CR(F(p),c\tau).
$$
\item Fixons $\lambda>0$ et $0<\eta<\lambda$. Il existe une constante 
$c(\eta)>0$ telle que :
$$
\begin{array}{l}
\hspace{-1cm}\forall\eps\leq\eps_0,\;\; \forall F\in \F_{p,a} \;:\\
F(A_\eps(p,\tau))\subset \{\Lambda_a\geq \lambda\} 
\Longrightarrow F(B^\CR(p,\tau))\supset B_W^\CR\big(a,\frac{c(\eta) 
\tau}{\sqrt\eps}\big)\cap \{\Lambda_a\leq \lambda-\eta\}.
\end{array}
$$
\item Soit $y\in W\priv\{a\}$ et $\eta<d(a,y)$. Il existe une constante 
$c(\eta)>0$  telle que :
$$
\begin{array}{l}
\hspace{-1cm}\forall\eps\leq\eps_0,\;\; \forall F\in \F_{p,a} \;:\\
F(A_\eps(p,\tau))\subset B(y,\eta)\Longrightarrow
F(B^\CR(p,\tau))\supset B_W^\CR\big(a,\frac{c(\eta)\tau}{\sqrt\eps}\big) 
\priv B(y,2\eta).
\end{array}
$$
\end{enumerate}\vspace{-0.1cm}
De plus, les constantes $c(K)$ et $c(\eta)$  ne 
dépendent que des valeurs de la forme de Lévi sur $U$ et $W$. Les points 1, 2 et 
3 restent vrais pour des constantes moitié lorsque $\om,\wdt \om$  sont remplacés 
par de petites perturbations au sens de Hausdorff pourvu que $U$ et $W$ soient 
eux faiblement perturbés au sens $\cc^2$.
\end{prop}
\begin{figure}[h]
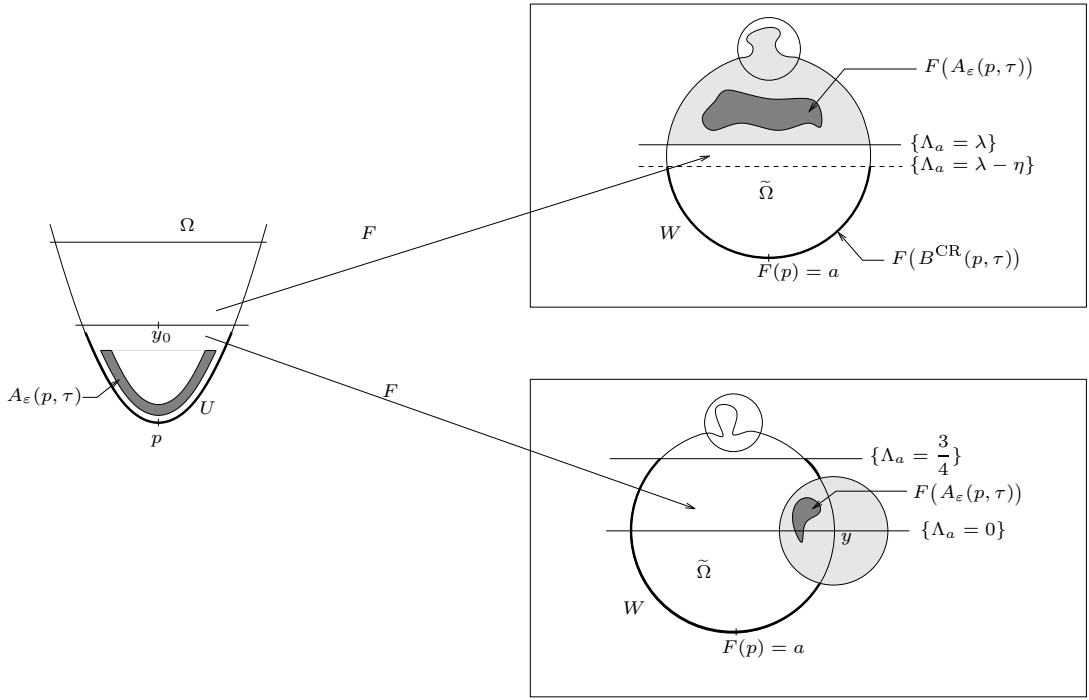

\begin{center}
\input popurri.pstex_t
\end{center}
\caption{Situation 2 (en haut), situation 3 (en bas)}
\label{popourri}
\end{figure}
\section{Démonstration du théorème 1.}
En  utilisant de façon systématique la proposition \ref{recap}, nous allons 
vérifier les propriétés $\K(F_n)$, $O(F_n)$, $\K(G_n)$ et $O(G_n)$. Ceci 
établira le théorème 1 {\it via} la proposition \ref{reduction}.
Toute la preuve du théorème 1 repose en fait sur  l'injectivité de 
$F_n$ au bord et les choix effectués lors de la renormalisation.
 
Rappelons que $\wdt\om_n$ converge vers $B$ au sens $\mathcal C^2$ sur les compacts 
de $\C^2\priv \{(1,0)\}$ (propriété \ref{cvcool}). Ainsi, pour $n\geq n_0$ 
suffisamment grand, $[b\wdt \om_n]_{3/4}$ est 
strictement pseudoconvexe et $\Lambda_x$ est positive sur $\wdt\om_n$ pour $x\in 
[b\wdt\om_n]_{3/4}$. Le domaine $\om$ étant quant à lui fixé et strictement convexe, 
la proposition \ref{recap} s'appliquera donc avec 
$\{U,W\}:=\{[b\om]_1,[b\wdt \om_n]_{3/4}\}$.
\subsection{Preuve de ``$\K \Longrightarrow O$''.}
\noindent\underline{$\K(F_n)\Longrightarrow O(F_n)$ :}
Appliquons le point 1 de la proposition \ref{recap} à $F_n$ en 
prenant $U:=b\om_1$ et $W_n:=[b\wdt\om_n]_{3/4}$. Comme $F_n(p)=a$ et 
qu'en vertu de l'hypothèse $\K(F_n)$ $F_n(A_{\eps_0}(p,\tau))$ tend 
vers $F(A_{\eps_0}(p,\tau))\Subset B$, on obtient :
$$
F_n\big(B^\CR(p,\tau)\big)\supset B_{W_n}^\CR(a,c\tau).
$$
Compte tenu du lemme \ref{remplissage} il suffit alors de vérifier qu'il existe 
$\lambda>0$ tel que $B^\CR_{W_n}(a,c\tau)\supset [b\wdt \om_n]_\lambda$ pour 
tout $n\in \N$. Or ceci résulte immédiatement du fait que $W_n$ converge au sens 
$\cc^2$ vers $[bB]_{3/4}$.\cqfd

Lorsque $O(F_n)$ est vérifiée, le lemme \ref{remplissage} assure l'existence 
d'un inverse à droite $G_n$ pour $F_n$ défini sur $\wdt \om_{n,\lambda}$. Cette 
suite d'applications est normale car elle prend ses valeurs dans $\om_1$, on peut 
donc la supposer convergente vers une application $G:B_\lambda\lra \overline{\om}_1$. 
Comme ci-dessus, nous allons montrer que si $G$ est à valeurs dans $\om_2$ alors 
$O(G_n)$ est satisfaite.\vspace{0.1cm}\\
\underline{$\K(G_n)\Longrightarrow O(G_n)$ :} 
Appliquons le point 1 de la proposition \ref{recap} à $G_n$ en 
prenant $U_n:=[b\wdt\om_n]_{\lambda}$ et $W:=[b\om]_2$. Comme $G_n(p)=a$ et 
qu'en vertu de l'hypothèse $\K(G_n)$ $G_n(A_{\eps_0}(a,\tau))$ tend 
vers $G(A_{\eps_0}(p,\tau))\Subset \om$, on obtient :
\begin{equation}\label{kog}
G_n\big(B_{U_n}^\CR(a,t)\big)\supset B_W^\CR(p,ct)\hspace{0.5cm} \forall 
t\leq\tau.
\end{equation}
Par ailleurs, la distance $\text{CR}$ majorant la distance euclidienne, 
on a $B_{U_n}^{\text{\tiny CR}}(a, t)\subset [b\wdt\om_{n}]_t$.
De plus, $B^{\text{\tiny CR}}(p,t)\supset [b\om]_{At^2}$ où $A$
désigne le minimum des valeurs propres de la forme de Levi de $[b\om]_1$ 
(il s'agit d'un cas particulier très simple du théorème 4 de \cite{nsw}). 
La propriété $O(G_n)$ découle alors de (\ref{kog}) et du lemme \ref{remplissage}.\cqfd

En remplaçant $p$ par un point arbitraire $q\in [b\om]_1$ dans la preuve de 
$\K(F_n)\Longrightarrow O(F_n)$, on obtient une propriété légèrement plus 
forte que $O(F_n)$. Celle-ci nous sera utile pour établir $\K(G_n)$.
\begin{lemme}\label{ptab}
Si $\K(F_n)$ est vraie alors il existe des constantes $c,\tau>0$ telles que 
pour tout $q\in [b\om]_1$ et $n$ suffisamment grand, 
$$
F_n(q)\in [b\wdt \om_n]_{3/4} \Longrightarrow F_n\big(B^\CR(q,\tau)\big)
\supset B_{[b\wdt \om_n]_{3/4}}^\CR(F_n(q),c\tau).
$$
\end{lemme}
\subsection{Preuve de $\K(F_n)$.}
En raison de l'alternative régissant la convergence de $F_n$ 
(corollaire \ref{normal}), 
il suffit d'établir que dans la situation décrite dans la première partie, 
$\wdt y_n=F_n(y_0)$ reste dans un compact de $B$.
Raisonnons par l'absurde, supposons que $\wdt y_n$ tende vers
 $y_\infty\in bB$. La suite  $(F_n)$ converge alors 
localement uniformément vers $y_\infty$ (corollaire \ref{normal}). 
A $\eps>0$ correspond donc un entier $n_\eps$ tel que 
$F_{n_\eps}(A_{\eps}(p,\tau)\big)\subset B(y_\infty,1/4)$. En appliquant 
le point 3 de la proposition \ref{recap} à $F_{n_\eps}$ avec $U:=[b\om]_1$,
$W_{n_\eps}:=[b\wdt\om_{n_\eps}]_{3/4}$, $\eta=1/4$ et $(a,y):=(a,y_\infty)$, 
on obtient :
\begin{equation}\label{kf1}
F_{n_\eps}([b\om]_1)\supset  
B^\CR_{W_{n_\eps}}\left(a,\frac{c(\eta)\tau}{\sqrt \eps}\right)
\,\priv B(y_\infty,\frac{1}{2}).
\end{equation}
Comme les $[b\wdt \om_n]_{3/4}$ tendent vers $[bB]_{3/4}$ au sens $\cc^2$, leurs 
diamètres $\text{CR}$ sont uniformément borné. Il s'ensuit que  pour $\eps$ 
suffisamment petit, 
(\ref{kf1}) donne :
\begin{equation} \label{kf2}
F_{n_\eps}([b\om]_1)\supset  [b\wdt\om_{n_\eps}]_{3/4}\,\priv B(y_\infty,\frac{1}{2}).
\end{equation}
Dorénavant, $\eps$ est fixé et nous noterons $n$ au lieu de $n_\eps$.
L'application ${F_n}_{|[b\om]_2}$ étant un difféomorphisme sur son image, 
$F_n([b\om]_1)$ est contractile. Alors (\ref{kf2}) montre que $F_n([b\om]_1)$ 
doit contenir $[b\wdt \om_n]_{3/4}$. Il vient donc :
\begin{equation}\label{eqlem1}
b\big[F_n([b\om]_1)\big]=F_n(b\om\cap\{\re z_1=1\})\subset [b\wdt \om_n]_{3/4}^+.
\end{equation} 
Cette inclusion conduit à une contradiction car le principe du maximum 
appliqué à la fonction  $-\re z_1\circ F_n$ montre que 
$F_n(\om\cap\{\re z_1=1\})\subset \wdt\om_{n,3/4}^+$ alors que par 
construction $F_n(y_0)\in \{\re z_1=0\}$. \hfill $\square$
\subsection{Preuve de $\K(G_n)$.}
Définissons $\lambda_n:=\sup\{t\in\R \;|\;F_n(\om_1)\supset\wdt\om_{n,t}\}$. 
Comme $F$ est un difféomorphisme local sur $\om_1$ (lemme \ref{remplissage}), 
le bord de $F_n(\om_1)$ dans $\wdt \om_n$ est inclus dans 
$F_n(\om \cap\{\re z_1=1\})$. Par définition de $\lambda_n$,
$$
\inf\big\{\re z_1\circ F_n(x),\;x\in \om\cap\{\re z_1=1\}\big\}=\lambda_n.
$$
 Le principe du maximum appliqué à $-\re z_1\,
\circ F_n$ sur $\om\cap\{\re z_1=1\}$ montre donc :
\begin{equation}\label{z'n}
\exists z'_n\in \{\re z_1=1\}\cap b\om, \hspace{0.5cm} 
F_n(z_n')\in\{\re z_1\leq\lambda_n\}.
\end{equation}
Puisque $\K(F_n)$ et donc $O(F_n)$ sont satisfaites, les valeurs d'adhérence de  
$\lambda_n$ sont supérieures à  $\lambda>-1$. Quitte à extraire, la suite des 
applications $G_n$ converge donc vers $G:B_{\lambda_\infty}\lra \overline{\om_1}$,
 où $\lambda_\infty$ désigne une valeur d'adhérence de $\lambda_n\geq \lambda$.

Lorsque $\lambda_\infty>0$, la suite $\wdt y_n$ est compacte dans 
$B_{\lambda_\infty}$ d'après $\K(F_n)$. De plus $G_n(\wdt y_n)=y_0\in \om_2$. 
L'application limite $G$ est donc à valeurs dans $\om_2$ : $\K(G_n)$ est vérifiée. 

Lorsque $\lambda_\infty\leq 0$, nous prouvons ci-dessous qu'il existe une suite 
$(x_n)$ relativement compacte dans $\om_1$ et une constante $\eta_0>0$ telles que 
$F_n(x_n)\in \wdt\om_{n,\lambda_\infty-\eta_0}$ pour $n$ assez grand. 
Alors $\{F_n(x_n)\}\Subset B_{\lambda_\infty}$ d'après $\K(F_n)$. On conclut
cette fois en utilisant $\big\{G_n [F_n(x_n)]\big\}=\big\{x_n\big\}\Subset \om_2$. 

\paragraph{} 
Fixons  $\tau>0$ tel que $B^{\text{\tiny CR}}(p,\tau)\Subset [b\om]_1$. Il suffit 
pour établir l'existence de $(x_n)$ de voir qu'il existe $\eps_0,\eta_0>0$ tels que :
$$
F_n(A_{\eps_0}(p,\tau))\cap \wdt \om_{n,\lambda_n-\eta_0}\neq \emptyset \hspace{0.5cm}
\forall n\gg 1.
$$
Raisonnons par l'absurde : supposons que pour tout  $\eps,\eta>0$, 
il existe un entier $n_\eps$ tel 
que $F_{n_\eps}(A_\eps(p,\tau))\subset \wdt\om_{n_\eps,\lambda_{n_\eps}-\eta}^+$. 
En appliquant le point 2 de la proposition \ref{recap}  à $F_{n_\eps}$ avec 
$U:=[b\om]_1$, $W:=[b\wdt\om_n]_{1/2}$ et $(\lambda,\eta):=(\lambda_{n_\eps},\eta)$
on obtient :
\begin{equation}\label{kg1}
F_{n_\eps}\big(B^\CR(p,\tau/2)\big)\supset  B^\CR_{W_{n_\eps}}
\left(a,\frac{c(\eta)\tau}{2\sqrt \eps}\right)\cap\{\re z_1\leq \lambda_{n_\eps}-\eta\}.
\end{equation}
Comme les $[b\wdt \om_n]_{3/4}$ tendent vers $[bB]_{3/4}$ au sens $\cc^2$, leurs 
diamètres $\text{CR}$ sont uniformément borné. Il s'ensuit que  pour $\eps$ 
suffisamment petit, (\ref{kg1}) donne :
$$
F_{n_\eps}\big(B^\CR(p,\tau/2)\big)\supset[b\wdt\om_{n_\eps}]_{\lambda_{n_\eps}-\eta}.
$$ 
Ceci étant valable pour $\eta$ arbitrairement petit, on obtient après extraction 
éventuelle :
\begin{equation}\label{fait1}
d^{\text{\tiny CR}}\left(bF_n\big(B^{\text{\tiny
CR}}(p,\tau/2)\big),[b\wdt\om_{n}]_{\lambda_n}^+\right)\underset{n\to\infty}{\lra}
0.  
\end{equation}
De plus, on voit à l'aide du lemme \ref{ptab} que la distance $\text{CR}$ entre les 
points de $bF_n(B^{\text{\tiny CR}}(p,\tau))$ et les points de 
$F_n(B^\CR(p,\tau/2))$ est bornée inférieurement par $c\tau/2$. 
Il vient donc de (\ref{fait1}) 
que $F_n(B^{\text{\tiny CR}}(p,\tau))$ contient $[b\wdt\om_{n}]_{\lambda_n}$ 
pour $n$ assez grand. Avec (\ref{z'n}), ceci contredit l'injectivité de 
${F_n}_{|U}$. \hfill $\square$ 
\section{Applications.}
Cette partie comporte trois applications du théorème 1. La première est un
résultat de Wong-Rosay local, la deuxième concerne l'étude de la dynamique des
applications holomorphes propres et la dernière un 
résultat bien connu de la géométrie $\text{CR}$. 
\subsection{Une version locale du théorème de Wong-Rosay.}
\begin{thm2}
Soient $\om$, $\om'$ deux domaines de $\C^k$ pseudoconvexes à bords lisses et $U$
un ouvert de $\spc(b\om)$. Soit $(f_n)_n$ une suite d'injections holomorphes de
$\om$ dans $\om'$ propres relativement à $U$. S'il existe un point $y_0\in \om$
dont l'orbite $(f_n(y_0))_n$ accumule un point $y\in \spc(b\om')$, alors $U$
est localement sphérique.
\end{thm2} 
\begin{rque}
Si les $f_n$ sont des surjections alors il s'agit du théorème classique.
\end{rque}
Par souci de lisibilité, nous appellerons dans cette preuve ordonnée d'un point
$z=(z_1,z')\in \C^k$ la partie réelle de $z_1$. \\
\noindent\underline{Preuve :} On peut bien entendu supposer que $U$ est connexe. 
Si $f_n$ converge localement uniformément sur $U$
vers $y$ alors le théorème 1 s'applique et assure que $U$ est sphérique. Sinon il 
existe $p\in U$ et des points $p_n\in U$ qui tendent vers $p$ tels que 
les $f_n(p_n)$ restent loins de $y$ (quitte à extraire). La sphéricité de $U$ découlera 
alors immédiatement des assertions suivantes :
\begin{eqnarray}
&\text{ la composante connexe $C_y$ de $y$ dans $\spc(b\om')$ est sphérique,}
\label{1'2} &\\
& f_n:U\lra C_y \text{ est un difféomorphisme } \text{CR}.\label{1'1}& 
\end{eqnarray}

Une estimation classique de la métrique
de Kobayashi au voisinage de $y$ \cite{graham} montre que $f_n$ 
converge uniformément vers $y$ sur les compacts de $\om$. Quitte à restreindre 
$\om$ et à effectuer un changement de variables, on peut donc supposer que $p=y=0$,
$T_pb\om=T_yb\om'=\{\re z_1=0\}$, $U=b\om_1$, $y_0=(1/2,0)$ et
$f_n(y_0)\in \om'_1$. On notera aussi $\eps_n$ l'ordonnée du point $p_n$.

Comme $y\in \spc(b\om')$, il existe une fonction \psh $\chi$ définie sur
$\om'$, continue sur $\overline{\om}'$ et présentant un maximum global strict
en $y$. Le principe du maximum appliqué à $\chi\circ f_n$  sur 
$\om_1\cap\{\re z_1=\eps_n\}$ montre qu'il existe des points $p'_n$
de mêmes ordonnées que $p_n$ tels que $f_n(p'_n)$ tend vers $y$. Or un 
argument dû à Fornaess \cite{fornaess} assure que les restrictions des $f_n$ à 
$U$ sont des difféomorphismes $\text{CR}$ et que $f_n(U)\subset \spc(b\om')$.
Ceci établit donc (\ref{1'1}). De plus, Quitte à modifier $p_n$ en le 
déplaçant le long d'un chemin d'ordonnée constante le joignant à $p'_{n}$, 
on peut supposer que $f_n(p_n)\in b\om'\cap\{\re z_1=\lambda\}$ et converge 
vers un point $a$. 

Pour prouver que $C_y$ est sphérique, il convient de noter que les fonctions
pics de classe $\cc^1$ que possèdent les points de $C_y$ peuvent remplacer 
les fonctions $\Lambda_q$ dans le lemme \ref{recouvre} et donc dans le point 3
de la proposition \ref{recap}. Un argument identique à celui utilisé pour 
prouver $\K(F_n)$ avec $(p_n,a,y_0,y)$ en place de 
$(p,a,y_0,y_\infty)$ montre que $f_n(b\om_\delta)$ contient 
$C_{y,\eps}:=\{p\in C_y\;|\; d(p,bC_y)>\eps\}$ pour $\delta,\eps$ arbitrairement
petits pourvu que $n$ soit assez grand. 
Ceci signifie que la suite des inverses $g_n$ des difféomorphismes $\text{CR}$
${f_n}_{|b\om_1}$  converge uniformément vers $p$ sur $C_{y,\eps}$. D'après le
théorème 1, $C_{y,\eps}$ est sphérique. On en déduit la sphéricité de $C_y$, 
c'est-à-dire (\ref{1'2}), en faisant tendre $\eps$ vers $0$.\hfill $\square$

\subsection{Dynamique des applications holomorphes propres.}
\begin{prop}\label{dynamique1}
Soit $\om$ un domaine convexe borné de $\C^k$ à bord lisse 
et $f$ une auto-application holomorphe propre de $\om$. S'il
existe un point $y_0\in \om$ dont l'orbite $\{f^n(y_0)\}$ accumule un point
$y$ de l'ensemble de stricte convexité $\SC(b\om)$ alors $\SC(b\om)$ est sphérique.
\end{prop}
\noindent \underline{Preuve :} Les théorèmes classiques de prolongement
au bord des applications holomorphes propres s'appliquent pour les domaines 
convexes à bords lisses : $f$ se prolonge en une application $\text{CR}$  de 
$b\om$ dans lui-même \cite{boasstraube} et $f:\SC(b\om)\to \SC(b\om)$
est un difféomorphisme local. De plus, il a été prouvé dans \cite{moi} 
que $f^n$ converge localement uniformément vers $y$ sur $\SC(b\om)$. Il suffit
donc de voir que tout point $p\in \SC(b\om)$ possède un voisinage sur lequel 
les $f^n$ sont toutes injectives pour que le théorème 1 s'applique et donne
le résultat souhaité.

Soit $V$ un voisinage de $y$ dans $C_y$
sur lequel $f$ est injective et $U_0$ un voisinage de $p$ relativement compact
dans $C_y$. Il existe un entier $n_0$ tel que $f^n(U_0)\subset
V$ pour tout $n\geq n_0$. Soit $U$ un voisinage de $p$ dans $U_0$
sur lequel  $f,f^2,\dots,f^{n_0}$ sont injectives. Montrons par récurrence 
que $f^n$ est injective sur $U$ pour tout $n$. 
Soient $x,x'\in U$ tels que $f^n(x)=f^n(x')$ avec
$n>n_0$. Alors $f(f^{n-1}(x))=f(f^{n-1}(x'))$ et $f^{n-1}(x),f^{n-1}(x')\in
V$. Comme $f$ est injective sur $V$, on en déduit que
$f^{n-1}(x)=f^{n-1}(x')$. En itérant ce procédé on obtient
$f^{n_0}(x)=f^{n_0}(x')$ donc $x=x'$.\cqfd
La proposition ci-dessus est valable pour d'autres types de domaines, dans lesquels
on peut assurer que l'itération de $f$ contracte $\spc(b\om)$ vers $a$ \cite{moi}.
En développant les idées du paragraphe précédent, on peut établir un résultat 
plus général. La technicité de sa preuve nous incite toutefois à ne pas la 
présenter ici.
\begin{theoreme}\label{thmchiant}
Soit $\om$ un domaine pseudoconvexe de $\C^k$ à bord lisse et $f$ une 
auto-application holomorphe propre de $\om$. S'il existe un point $y_0\in \om$ 
dont l'orbite accumule un point de stricte pseudoconvexité alors $\spc(b\om)$ 
est localement sphérique. 
\end{theoreme}
\subsection{Etude métrique des hypersurfaces strictement pseudoconvexes
  non-sphériques.} 
Nous redémontrons ici un résultat bien connu de la géométrie $\text{CR}$
(voir \cite{webster, burshn,beloshapka,loboda}) concernant le groupe
d'isotropie d'un point $a$ d'une hypersurface strictement pseudo-convexe
non-sphérique de $\C^k$, c'est à dire le groupe des germes de difféomorphismes
$\text{CR}$ de $(M,a)$. Ce groupe occupe une place importante dans la classification
$\text{CR}$ des hypersurfaces puisqu'il joue un rôle de stabilisateur. Le
théorème ci-dessous a deux conséquences directes : les hypersurfaces
non-sphériques peuvent être pourvues d'une métrique $\text{CR}$-invariante
 et d'un champs de
vecteur $\text{CR}$-invariant transverse à la distribution complexe.
\begin{theoreme}
Soit $(M,a)$ un germe d'hypersurface strictement pseudoconvexe
non-sphérique. Soit $G$ le sous-groupe de GL$(T_a M)$ défini par :
$$
G:=\{f'(a), \text{ où }f:(M,a)\lra (M,a) \text{ est un germe de
  difféomorphisme $\text{CR}$}\} 
$$
Il existe un changement de base $\text{CR}$ de $T_a M$, 
\textit{i.e.} $P\in GL(T_a M)\cap GL_\C(T^\C_a M)$
% avec $P_{|T_a^\C M}\in \L_\C(T_a^\C M)$, 
tel que :
$$
G=P^{-1}\left\{
\left(
\begin{array}{cc}
1 & 0\\ 0 & A
\end{array}\right), \; A\in U(k-1)
\right\}P.
$$
{\it On dit que le groupe structural de la géométrie $\text{CR}$ des
  hypersurfaces strictement pseudoconvexe non-sphériques de $\C^k$ est
  réductible à  $U(k-1)$.}
\end{theoreme}
\noindent \underline{Preuve :}
Nous noterons ici $\L(a,u)$ la forme de Levi de $M$ en $a$ appliquée à un
vecteur $u\in T^\C_a b\om$ et $\kappa_-,\kappa_+$ ses valeurs
propres minimales et maximales (si $u$ désigne un champ de vecteur holomorphe
au voisinage de $a$ prolongeant $u$, $\L(a,u)=\langle [u,\bar u](a),i\vec
N(a)\rangle$). Soit $f$ un germe de difféomorphisme $\text{CR}$ de $(M,a)$. En
accord avec la partie 2, nous notons $n_f(a):=\langle f'(a)i\vec N(a),i\vec
N(a)\rangle$ ($n_f(a)$ est évidemment positif). Un calcul simple montre :
$$
\L(a,f'(a)u)=n_f(a)\L(a,u) \hspace{0.5cm} \forall u\in T^\C_a M.
$$ 

Si $n_f(a)<1$, il existe un entier $n_0$ tel que
$n_f(a)^{n_0}\kappa_+\kappa_-^{-1}<1$. Le germe de difféomorphisme
$f^{n_0}$ vérifie alors :
$$
\begin{array}{l}
\langle {f^{n_0}}'(a)i\vec N(a),i\vec N(a)\rangle=n_f(a)^{n_0}<1,\\
\Vert {f^{n_0}}'(a)u\Vert^2\leq \kappa_-^{-1} \L(a,{f^{n_0}}'(a)u)\leq
n_f(a)^{n_0}\kappa_+\kappa_-^{-1} \Vert u\Vert^2<\Vert u\Vert^2,
\end{array}
$$  
si bien que $f^{n_0}$ est un germe de difféomorphisme contractant de
$(M,a)$. L'hypothèse de non-sphéricité de $M$ assure {\it via}  le théorème 1
de l'impossibilité de cette situation ainsi évidemment de celle où
$n_f(a)>1$. On conclut donc que $n_f(a)=1$. Ceci a pour conséquence que
$\L(a,f'(a)u)=\L(a,u)$ pour $u\in T^\C_aM$ donc $f'(a)_{|T^\C_aM}$ est une
isométrie de la forme hermitienne $\L(a,\cdot)$. Si $\bc^\C$ désigne une base
orthonormée pour $\L(a,\cdot)$, et $\bc=(i\vec N(a),\bc^\C)$, on a donc :
$$
\text{Mat}(f'(a),\bc,\bc)=
\left(
\begin{array}{cc}
1 & 0\\
* & A
\end{array}
\right) \hspace{0.5cm} \text{où } A\in U(k-1).
$$
Il suffit alors de trouver un vecteur invariant de $i\vec N+T^\C_a M$ sous
l'action de $G$ pour que la matrice de $f'(a)$ dans la base constituée de ce
vecteur et de $\bc^\C$ soit telle que souhaitée. 
 
A cette fin, remarquons que $G$ induit une action affine sur l'espace 
$\ec:=\{i\vec N+T^\C_a M\}$ qui stabilise le compact $K:=\{[u,\bar u](a), u\in
T^\C_a M, \L(a,u)=1\}$ ($K\subset \ec$ par définition de la forme de Levi). Le
 centre d'inertie de $K$ est donc un point fixe sous l'action de $G$ du type 
souhaité.\hfill $\square$ 

{
\footnotesize
\bibliography{bib2.bib}

\begin{thebibliography}{10}

\bibitem{bedpin1}
E.~Bedford and S.~Pinchuk.
\newblock Domains in {${\bf C}\sp {n+1}$} with noncompact automorphism group.
\newblock {\em J. Geom. Anal.}, 1(3):165--191, 1991.

\bibitem{bedpin2}
E.~Bedford and S.~I. Pinchuk.
\newblock Domains in {${\bf C}\sp 2$} with noncompact groups of holomorphic
  automorphisms.
\newblock {\em Mat. Sb. (N.S.)}, 135(177)(2):147--157, 271, 1988.

\bibitem{bedpin3}
E.~Bedford and S.~I. Pinchuk.
\newblock Convex domains with noncompact groups of automorphisms.
\newblock {\em Mat. Sb.}, 185(5):3--26, 1994.

\bibitem{bell}
S.~Bell.
\newblock Local boundary behavior of proper holomorphic mappings.
\newblock In {\em Complex analysis of several variables (Madison, Wis., 1982)},
  volume~41 of {\em Proc. Sympos. Pure Math.}, pages 1--7. Amer. Math. Soc.,
  Providence, RI, 1984.

\bibitem{beloshapka}
V.~K. Beloshapka.
\newblock The dimension of the group of automorphisms of an analytic
  hypersurface.
\newblock {\em Izv. Akad. Nauk SSSR Ser. Mat.}, 43(2):243--266, 479, 1979.

\bibitem{berteloot2}
F.~Berteloot.
\newblock Characterization of models in {$\bold C\sp 2$} by their automorphism
  groups.
\newblock {\em Internat. J. Math.}, 5(5):619--634, 1994.

\bibitem{berteloot3}
F.~Berteloot.
\newblock Attraction des disques analytiques et continuit\'e h\"old\'erienne
  d'applications holomorphes propres.
\newblock In {\em Topics in complex analysis}, volume~31 of {\em Banach Center
  Publ.}, pages 91--98. 1995.

\bibitem{berteloot}
F.~Berteloot.
\newblock Méthodes de changement d'échelles en analyse complexe.
\newblock {\em Ann. Fac. Sci. Toulouse Math.}, 2005.

\bibitem{boasstraube}
H.~Boas and E.~Straube.
\newblock Sobolev estimates for the {$\overline\partial$}-{N}eumann operator on
  domains in {${\bf C}\sp n$} admitting a defining function that is
  plurisubharmonic on the boundary.
\newblock {\em Math. Z.}, 206(1):81--88, 1991.

\bibitem{bogess}
A.~Boggess.
\newblock {\em C{R} manifolds and the tangential {C}auchy-{R}iemann complex}.
\newblock Studies in Advanced Mathematics. CRC Press, Boca Raton, FL, 1991.

\bibitem{burshn}
D.~Burns, Jr. and S.~Shnider.
\newblock Geometry of hypersurfaces and mapping theorems in {${\bf C}\sp{n}$}.
\newblock {\em Comment. Math. Helv.}, 54(2):199--217, 1979.

\bibitem{chemos}
S.~S. Chern and J.~K. Moser.
\newblock Real hypersurfaces in complex manifolds.
\newblock {\em Acta Math.}, 133:219--271, 1974.

\bibitem{fornaess}
J.~E. Fornaess.
\newblock Biholomorphic mappings between weakly pseudoconvex domains.
\newblock {\em Pacific J. Math.}, 74(1):63--65, 1978.

\bibitem{gaussier}
H.~Gaussier.
\newblock Characterization of models for convex domains.
\newblock {\em Preprint}.

\bibitem{graham}
I.~Graham.
\newblock Boundary behavior of the {C}arath\'eodory and {K}obayashi metrics on
  strongly pseudoconvex domains in {$C\sp{n}$} with smooth boundary.
\newblock {\em Trans. Amer. Math. Soc.}, 207:219--240, 1975.

\bibitem{loboda}
A.~V. Loboda.
\newblock Linearizability of automorphisms of nonspherical surfaces.
\newblock {\em Izv. Akad. Nauk SSSR Ser. Mat.}, 46(4):864--880, 1982.

\bibitem{nsw}
A.~Nagel, E.~M. Stein, and S.~Wainger.
\newblock Balls and metrics defined by vector fields. {I}. {B}asic properties.
\newblock {\em Acta Math.}, 155(1-2):103--147, 1985.

\bibitem{moi}
E.~Opshtein.
\newblock Dynamique des applications holomorphes propres des domaines
  r{\'e}guliers et probl{\`e}me de l'injectivit{\'e}.
\newblock {\em Submitted}, 2004.

\bibitem{pinchuk}
S.~Pinchuk.
\newblock The scaling method and holomorphic mappings.
\newblock In {\em Several complex variables and complex geometry, Part 1},
  volume~52 of {\em Proc. Sympos. Pure Math.}, pages 151--161. Amer. Math.
  Soc., 1991.

\bibitem{pinchuk3}
S.~I. Pinchuk.
\newblock Proper holomorphic maps of strictly pseudoconvex domains.
\newblock {\em Sibirsk. Mat. \v Z.}, 15:909--917, 959, 1974.

\bibitem{pinchuk4}
S.~I. Pin{\v{c}}uk.
\newblock Proper holomorphic mappings of strictly pseudoconvex domains.
\newblock {\em Dokl. Akad. Nauk SSSR}, 241(1):30--33, 1978.

\bibitem{rosay}
J.-P. Rosay.
\newblock Sur une caract\'erisation de la boule parmi les domaines de {${\bf
  C}\sp{n}$} par son groupe d'automorphismes.
\newblock {\em Ann. Inst. Fourier (Grenoble)}, 29(4):ix, 91--97, 1979.

\bibitem{webster}
S.~M. Webster.
\newblock On the transformation group of a real hypersurface.
\newblock {\em Trans. Amer. Math. Soc.}, 231(1):179--190, 1977.

\bibitem{wong}
B.~Wong.
\newblock Characterization of the unit ball in {${\bf C}\sp{n}$} by its
  automorphism group.
\newblock {\em Invent. Math.}, 41(3):253--257, 1977.

\end{thebibliography}
\bibliographystyle{abbrv}
}
\end{document}